\documentclass[3p,times]{elsarticle}

\usepackage{amssymb}
\usepackage[figuresright]{rotating}

\usepackage{color}

\usepackage{pgfplots}
\pgfplotsset{compat=1.15}
\usetikzlibrary{arrows}

\def\RR{\hbox{\sf I\kern-.14em\hbox{R}}}

\newtheorem{Theorem}{Theorem} 
\newtheorem{lem}[Theorem]{Lemma}
\newtheorem{corro}[Theorem]{Corollary}

\newtheorem{proposition}[Theorem]{Proposition}

\newcommand{\be}{\begin{equation}}
\newcommand{\ee}{\end{equation}}

\newcommand{\DA}[1]{   {\color{black} #1}}
\newcommand{\JZ}[1]{   {\color{black} #1}} %{\color{blue} #1}}

\begin{document}

\begin{frontmatter}

%% Title, authors and addresses

%% use the tnoteref command within \title for footnotes;
%% use the tnotetext command for the associated footnote;
%% use the fnref command within \author or \address for footnotes;
%% use the fntext command for the associated footnote;
%% use the corref command within \author for corresponding author footnotes;
%% use the cortext command for the associated footnote;
%% use the ead command for the email address,
%% and the form \ead[url] for the home page:
%%
%% \title{Title\tnoteref{label1}}
%% \tnotetext[label1]{}
%% \author{Name\corref{cor1}\fnref{label2}}
%% \ead{email address}
%% \ead[url]{home page}
%% \fntext[label2]{}
%% \cortext[cor1]{}
%% \address{Address\fnref{label3}}
%% \fntext[label3]{}

%\dochead{}
%% Use \dochead if there is an article header, e.g. \dochead{Short communication}

%% use optional labels to link authors explicitly to addresses:
%% \author[label1,label2]{<author name>}
%% \address[label1]{<address>}
%% \address[label2]{<address>}

%\begin{center}
%{\bf Janez \v{Z}erovnik} \\
%FS, University of Ljubljana, A\v{s}ker\v{c}eva 6, 1000 Ljubljana, Slovenia\\
%%Institute of Mathematics, Physics and Mechanics, Jadranska 19, Ljubljana, Slovenia\\
%            and \\ 
%\Rev{Rudolfovo - scientific and technological center Novo mesto, Podbreznik 15, 8000 Novo mesto, Slovenia}\\
%            janez.zerovnik@fs.uni-lj.si \\ 
%}
%\end{center}

\title{On rainbow domination of cubic graphs}
\author{Janez \v{Z}erovnik}

\address{FS, University of Ljubljana, A\v{s}ker\v{c}eva 6, 1000 Ljubljana, Slovenia\\
%Institute of Mathematics, Physics and Mechanics, Jadranska 19, Ljubljana, Slovenia\\
            and \\ 
Rudolfovo - scientific and technological center Novo mesto, Podbreznik 15, 8000 Novo mesto, Slovenia\\
            janez.zerovnik@fs.uni-lj.si \\ 
}

\begin{abstract}
The structure of minimal weight rainbow domination  functions of cubic graphs are studied. 
Based on general observations for cubic graphs, generalized Petersen graphs $P(ck,k)$ are characterized whose 
4- and 5-rainbow domination numbers equal the general lower bounds.  
As $t$-rainbow domination of cubic graphs  for $t \ge 6$ is trivilal, 
characterizations of such  generalized Petersen graphs $P(ck,k)$  are known for all $t$-rainbow domination numbers.
In addition, new upper bounds  for 4- and 5-rainbow domination numbers that are valid for all $P(ck,k)$ are provided. 
\end{abstract}

\begin{keyword}
rainbow domination; rainbow  domination number;  cubic graphs; generalized Petersen graphs
%% keywords here, in the form: keyword \sep keyword

%% MSC codes here, in the form: \MSC code \sep code
%% or \MSC[2008] code \sep code (2000 is the default) 

\end{keyword}

\end{frontmatter}

%%
%% Start line numbering here if you want
%%
% \linenumbers

%% main text

%% The Appendices part is started with the command \appendix;
%% appendix sections are then done as normal sections
%% \appendix

%% \section{}
%% \label{}

%% References
%%
%% Following citation commands can be used in the body text:
%% Usage of \cite is as follows:
%%   \cite{key}         ==>>  [#]
%%   \cite[chap. 2]{key} ==>> [#, chap. 2]

\section{Introduction}

 A simple graph $G = (V(G), E(G))$ is a combinatorial object, where 
 $V=V(G)$ is a set whose elements are called  a set of {\em vertices} and  with $E=E(G)$  a set of {\em edges}. 
 Edges in simple undirected graphs are pairs of vertices, $e= \{u,v\}\in E(G)$.   
% (We often write shortly $uv$ instead of $ \{u,v\}$.) 
In such case we say that vertices $u$ and $v$ are {\em  neighbors}.
 The set of all neighbors of a vertex is its  {\em neighborhood}. The number of its neighbors is called the {\em degree} of a vertex.
 A graph is called {\em 3-regular} or {\em cubic} if all vertices in $V(G)$ are of degree 3. 

A $t$-rainbow domination function ($t$RDF)  assigns   a subset
of the   set  of colors $\{1,2,\cdots,t\}$   to each  vertex of $G$ such that
 each vertex that is assigned  an empty set has all $t$ colors in its  neighborhood.
The weight of assignment  $g$ is  the value $w(g) =\sum_{v \in V(G)}    w(g(v))  $,   
where $w(g(v))$ is the number of colors assigned to vertex $v$.
We also say that $G$ is \emph{ $t$RD-colored} (or simply, colored) by $g$.   
A  vertex is said to be  $t$RD-dominated  if
either: 
(1) it is assigned a  nonempty set of colors,
or, 
(2) it has all colors in its neighborhood.
  If $g(v) \neq \emptyset$, a vertex $v$ is said to be \emph{colored},   and is not colored or \emph{uncolored} otherwise.
The minimum weight over all $tRD$ functions of $G$ is called the \emph{$t$-rainbow domination number} $\gamma_{rt}(G)$.
If the sets assigned by $g $ have cardinality at most one,     $g$ is a  {\em singleton  $tRD$ function}   (see \cite{ErvesSymmetry3}).
The minimal weight obtained when considering only  singleton  $tRD$ functions
 is called the {\em singleton $t$-rainbow domination number},   and is denoted by $\tilde\gamma_{rt}$.
Directly from definitions it follows, for any graph $G$ and any $t$,  $\gamma_{rt}(G) \le \tilde\gamma_{rt} (G)$. 

Rainbow domination has been first studied    in 2005   \cite{Bresar2005}  and has already received considerable popularity, perhaps because it si both 
an appealing concept and  motivated by practical applications. A survey \cite{Bresar2020} has appeared recently.
The invariant already has a number of varieties, for example rainbow independent domination \cite{KranerŠumenjak2018353}.

Rainbow domination of regular graphs is studied in \cite{Kuzman} 
where, among other results,  some general bounds are established for so called rainbow regular graphs. 
Here we generalize the bounds to cubic regular graphs that need not be rainbow regular.
Then use the general bounds to provide new bounds and in some cases exact values for rainbow domination of some families of generalized Petersen graphs.

Among the special classes of graphs that enjoy substantial interest in the community  are the generalized Petersen graphs 
\cite{Shao2019370,Gao2020a}. 
%\cite{Shao2019370,Gao2020a,                    .............................................. ADDrefs}.   ???
Besides the Petersen graphs $P(n,k)$ with small $k$,   the   Petersen graphs  $P(ck,k)$ received considerable attention. 
In \cite{Darja3RD},  bounds for 3-rainbow domination of the families $P(ck,k)$ were provided, 
generalizing the  results  of \cite{ErvesSymmetry3} for $P(6k,k)$ that include exact values for some subfamilies and tight bounds for other subfamilies. 
For analogous results regarding 2-rainbow domination of  $P(ck,k)$, see  \cite{ErvesSymmetry2,Simon2RD}.
Furthermore, note that 1-rainbow domination is equivalent to ordinary domination, which has been addressed in \cite{WangDOMINATION}.
Hence, for $r=1,2,3$ the $r$-rainbow domination of  $P(ck,k)$  is fairly well  studied  in the sense that we have 
exact values for some families and good lower and upper bounds for all other cases.
In this work, we first prove a  generalization of results of \cite{Kuzman}  from so called rainbow domination regular graphs to any regular graphs. 
Using this, we are able to extend the bounds and, in some cases, exact values to   $r$-rainbow domination for $r>3$. 

Main results in this paper are the following two theorems.
First theorem is a characterization of the generalized Petersen graphs that attain the general lower bound.

\begin{Theorem}    \label{MainTheoremX2novi} 
 Let   $t\in\{3,4,5\}$. 
Then $  \gamma_{rt}(P(n, k))   =   \frac{t}{3}n  $ if and only if  $n\equiv 0~\pmod 6$,   and $k \equiv 1, 5~\pmod 6$. 
%and, otherwise, 
%\begin{equation}
%  \gamma_{rt}(P(n, k))   >   \frac{t}{3} n \,.
%\end{equation} 
\end{Theorem} 
 
The second theorem provides lower and upper bounds for a subfamily of generalized Petersen graphs $P(ck,k)$.

\begin{Theorem}    \label{MainTheoremX1} 
Let  $t\in\{3,4,5\}$. Then 
\begin{equation}
 \frac{t}{3}ck   \le \gamma_{rt}(P(ck, k))    \le   \frac{t}{3} (c+1)(k+1)  - t \,.
\end{equation}
and 
$
  \gamma_{rt}(P(ck, k))   =   \frac{t}{3} ck   $  iff $c\equiv 0~\pmod 6$,   and $k \equiv 1, 5~\pmod 6$.
\end{Theorem} 
 
Theorem \ref{MainTheoremX1}   summarizes  the bounds of 
Theorems 3, 4, and 5  of  \cite{Darja3RD} and  provides a generalization  to 4 and 5-rainbow domination.
Furthermore,   Theorem   \ref{MainTheoremX2novi}  generalizes  the case in 
Theorem 3     \cite{Darja3RD} where exact value is obtained.
Using terminology of \cite{Kuzman},   Theorem   \ref{MainTheoremX2novi}  provides examples of 3-rainbow domination regular graphs. 
We just mention in passing that \cite{ZEROVNIK2024144} provides a characterization of generalized  Petersen graphs that are 
3-rainbow domination regular graphs.  
 Theorem  \ref{MainTheoremX2novi}  generalizes the observation  as it  implies that  the  neccessary and sufficient 
condition at the same time charaterizes the minimal examples for $t=3,4$ and 5. 

Theorem  \ref{MainTheoremX1}   extends  the  known  formulae and  bounds for  2-rainbow and 3-rainbow domination of $P(ck,k)$ \cite{Simon2RD,Darja3RD}
to $t =4$ and 5. For $t \ge 6$, the rainbow domination is well known \cite{Kuzman}  to be $|V(P(ck, k))|   = 2n$.

The rest of the paper is organized as follows.
In the next section, a short overview of related pevious work is given, emphasizing the results that are used later.
Section 3 provides observations that give some insight into the structure of rainbow dominating sets and lead to general lower and upper bounds for rainbow domination numbers of cubic graphs.
In Section 4  we study in more detail the generalized Petersen graphs. In particular, we prove the two Theorems given above. 
The last section give concluding remarks.

%%%%%%%%%%%%%%%%%%%%%%%%%%%%%%%%%%%%%%%%%%%%%%%%%%%%%%%%%%%%%%%%%%%%%%%%%%%%%%%%%
\section{Related Previous Work}

  Various  results on $k$-rainbow domination are already provided in the early  papers~\cite{Bresar2005, Bresar20072394, Bresar2008}.  
The problem is well-known to be NP-hard for general graphs \cite{Chang20108}.
In \cite{Pai} the authors provide an exact algorithm 
and a faster heuristic algorithm to  calculate  the 3-rainbow domination number. 
Therefore, in general, 3-rainbow domination numbers for small or moderate size graphs
can be computed, but it is  very hard or intractable  to handle large graphs.
Because of the  hardness of the general problem, 
it is interesting to study the complexity of the problem on restricted domains (c.f. trees)
and to consider particular graph classes. 
For example, it is known that the problem is NP-hard even when restricted to chordal graphs and to bipartite graphs, and there is a linear time algorithm for the $k$-domination problem on trees
\cite{Chang20108}. 
 
The special cases,  2-rainbow and 3-rainbow domination,   have been  studied  often in recent years.
In particular, the rainbow domination numbers $\gamma_{r2}$ and $\gamma_{r3}$ of several graph classes were established; 
see~\cite{Derya,wu2010bounds,WU2010706,Shao2019370,NatAcademy,Xu20092570} and the references~therein. 
In particular,  $k$-rainbow domination number of the Cartesian product of cycles, 
$C_n \square C_m$, for $k \geq 4$ is considered in \cite{Gao2019}.
Among other things, based on the results in  \cite{Amjadi2017394}, it is shown that 
$\gamma_{rk}(C_n \square C_m)=mn$ for $k \geq 8$. 
%V \cite{Gao1} je več natančnejših rezultatov za $4 \leq k \leq 7$.
%
In \cite{Gao2020c}, exact values of the 3-rainbow domination number of 
$C_3 \square C_m$ and $C_4 \square C_m$  
and   bounds on %the 3-rainbow domination number of 
$\gamma_{r3}(C_n \square C_m)$ for $n \geq 5$ are given.
In \cite{Fujita2015601}, sharp upper bounds on the $k$-rainbow domination number \DA{$\gamma_{rk}$} for all values of $k$
are proved. Even more,  the problem with minimum degree restrictions on the graph has been considered.
In particular, 
it was shown that for every connected graph $G$ of order $n \geq 5$, 
$\gamma_{r3}(G) \leq \frac{8n}{9}$.
In \cite{FURUYA201845}, the authors  prove that for every connected graph $G$ of order
$n \geq 8$ with minimal degree $\delta (G)\geq 2$, we have $\gamma_{r3}(G)\leq \frac{5n}{6}$.
 
In the past, generalized Petersen graphs have been studied extensively,
in many cases  as  counterexamples to 
conjectures  or as 
very interesting examples in research of various graph invariants.  
Often, subfamilies of generalized  Petersen graphs are considered.
Popular examples are graphs  $P(n,k)$ with fixed (and usually small) $k$, 
and 
$P(ck,k)$, for~fixed  $c$ and arbitrary $k$ (hence infinitely many $n=ck$).  
First bounds for 2-rainbow domination numbers of $P(n,3)$ have beed provided in \cite{Xu20092570}, 
namely it was shown that   $\gamma_{r2} (P(n,3)) \le  n - \lfloor \frac{n}{8} \rfloor + \beta$, 
where $\beta$ is either 0 or 1.
Later in \cite{Shao2019370}, it was proved that $\gamma_{r2}(P(n,1)) =n$ for $n\ge 5$. 
In the same paper,  exact values  for  $\gamma_{r3}(P(n,2))  $  and  upper bounds for  $\gamma_{r3}(P(n,3))  $ were obtained.
Exact values for   $\gamma_{r4}(P(n,1))  $  and  upper bounds for  $\gamma_{r5}(P(n,1))  $ are given in  \cite{Gao2020a}.
Furthermore, exact values of $\gamma_{rt}(P(n,1))$  for any $t \geq 8$ and $t = 4$ were derived and 
it is  proved that $\gamma_{rt}(P(2k,k))= 4 k$ for $t \geq 6$.  
The 3-rainbow domination numbers of some special classes of graphs such as paths, cycles and the generalized Petersen graphs $P(n,k)$ were investigated in \cite{Shao2014225}.
In particular,  3-rainbow domination numbers of $P(n,k)$ for some cases are established 
 and  the upper bounds for $P(n,2)$, $n\geq 5$, and $P(n,3)$,  $n\geq 30$, are provided.
Furthermore,  the general lower bound for 3-rainbow domination number  was established, 
$\gamma_{r3}(P(n,k))\geq n$, and
it was proved  that in case $k \equiv 1 \pmod 6$, $n \equiv 0 \pmod 6$ and $n>2k\geq 6$,  equality $\gamma_{r3}(P(n,k)) = n$  holds. 
In addition, it was determined that  for $n \geq 6$,
$\gamma_{r3}(P(n,1))= n +\alpha$, where 
$\alpha =0$ for $n \equiv 0 \pmod 6$,  
$\alpha =1$ for $n \equiv 1,2,3,5 \pmod 6$, and
$\alpha =2$ for $n \equiv 4 \pmod 6$.
The upper bound,
$\gamma_{r3}(P(n,2))\leq \lceil \frac{6n}{5}\rceil$ 
for $n \geq 5$   was proved.  
It follows % \cite{3RDNmain}, 
that 
$\gamma_{r3}(P(6k,k)) \geq 6k$ for each $k\geq 1$, 
$\gamma_{r3}(P(6k,k)) = 6k$ if $k \equiv 1 \pmod 6$,
 and 
$12\leq \gamma_{r3}(P(12,2)) \leq 15$.
 
Two-rainbow domination of generalized Petersen graphs $P(ck,k)$ was studied in \cite{Simon2RD}, 
providing a characterization of minimal examples and upper bounds for general case.
Three-rainbow domination of  $ P(ck,k) $  was studied in \cite{Darja3RD}, generalising the results on 3-rainbow domination of  $ P(6k,k) $  \cite{Erves2021b}.
The upper bounds for $\gamma_{r3}(P(6k,k))$ will be recalled in more detail later.
 
In \cite{Kuzman}, rainbow domination regular graphs are defined. 
More precisely,
a graph $G$ is $t$-rainbow domination regular ($t$-RDR graph) if it is $t$-regular and   $  \gamma_{rt} (G) = \frac{1}{2} |V(G)| $.
Let us recall three statements   regarding the rainbow domination of regular graphs  from \cite{Kuzman}  :

\begin{Theorem}  \label{LBKuzman}    \cite{Kuzman}
Let  $G$ be a $d$-regular graph. Then, for $\ell <2d$,  %we have 
%\begin{equation} 
$\displaystyle   \gamma_{r\ell}(G)   \ge  \left\lceil \frac{\ell  |V(G)| }{2d} \right\rceil
$
%\end{equation}
and 
$  \gamma_{r\ell}(G)   = V(G)$ for   $\ell \ge 2d$.
\end{Theorem} 
 
\begin{proposition}  \label{TheoremGeneral}      \cite{Kuzman}
Let  $G$ be a $d$-regular graph and  $d \le \ell <2d$. If   $  \gamma_{r\ell}(G)  = \frac{\ell |V(G)|}{2d}$, 
then $ \gamma_{r(\ell+1)}(G)   =    \frac{\ell+1 }{\ell}  \gamma_{r\ell}(G)   = \frac{(\ell+1) |V(G)|}{2d}$. 
\end{proposition} 

\begin{corro}  \label{PosledicaKuzman}      \cite{Kuzman}
Let  $G$ be a $d$-regular graph of order $n$. 
Then $G$ is $d$-rainbow regular if and only if 
$    \gamma_{r\ell}(G)   =   \frac{\ell  |V(G)| }{2d}  
$
for all  $d \le k \le 2d$.
\end{corro}

\section{Bounds for rainbow domination of cubic graphs}

 %%%%%%%%%%%%%%%%%%%%%%%%%%%%%%%%%%%%%%%%%%%%%%%%%%%%%%%%%%%%%%%%%%%%%%%%%%%%%%%%%
\subsection{Lower bounds for  3-regular graphs}
 
Let us introduce some more notation.
A 4RDF of $G$  $f$ with $n= |V(G)|$
gives rise to partition $V= V_0 \cup V_1 \cup V_2 \cup V_3  \cup V_4$ where $V_i$ denotes the set of vertices $v \in V_i$ with $w(f(v))= i$. 
Hence 
$n= n_0 + n_1 + n_2 + n_3 + n_4$ and $ w(f) =  n_1 + 2n_2 + 3n_3 + 4n_4$. 
In general, $n = \sum_{i=0}^{r} n_i= \sum_{i=0}^{r} |V_i|$ and  $w(f) = \sum_{i=0}^{r} n_i i$ for a $r$-rainbow domination function $f$.

\begin{lem}  \label{SplosnaMeja4lema} 
 Let $G$ be a 3-regular graph on $n$ vertices and $f$ be a $4RD$ function. Then  $  w(f) <  \frac{3}{4}|V(G)| $  implies  $n_3=n_4=0$. 
\end{lem}   

\begin{proof}
Assume that  the lemma  does not hold. 
Let $G$ be a minimal counterexample, i.e.  $  w(f) <  \frac{3}{4}|V(G)| $ for some 4RDF $f$   with $|V(G)| =n$ and   $n_3\neq 0$ or $n_4 \neq 0$. 
For the sake of contradiction, suppose that there is a vertex $v \in V(G)$ with $w(f(v)) = w_v \in\{3,4\}$.
Let us  construct a new graph $G'$ as follows.
First, delete the vertex $v$ and its neighbors. For each  neighbor $u$ of $v$, connect the two remaining neighbors with an edge. 
If $f(u) \neq  \emptyset$, move the colors of $u$ to any of its two neighbors.
By construction, $f$ restricted to $G'$ is a 4RDF, 
$w(f(G')) = w(f(G)) - w_v$, and $n' = |V(G')| =n -4$, so 
$$  w(f(G'))   = w(f(G)) - w_v < \frac{3}{4}n - w_v = \frac{3}{4} (n' +4) - w_v  =  \frac{3}{4} n' +  3 - w_v   \le   \frac{3}{4}n'  \, ,
$$
%and hence $G'$ is a counterexample with $n' <n$,
 contradicting the minimality of $G$.
\qed
\end{proof}  

\medskip\noindent
Lemma  \ref{SplosnaMeja4lema} gives some insight into the structure of minimal  weight 4RDR functions of arbitrary cubic graphs
that are not neccessarily rainbow domination regular.
The lemma also easily implies  the  next proposition that gives a lower bound for  $ \gamma_{r4} $. 

\begin{proposition}  \label{SplosnaMeja4} 
 Let $G$ be a 3-regular graph on $n$ vertices and $f$ be a $4RD$ function. Then $  w(f) \geq  \frac{2}{3}|V(G)|    \, .$ 
\end{proposition}   

\begin{proof}
Assume that $f$ is a 4RDF  of $G$ with  $n_3=n_4=0$. 
We are going to show that $w(f) \ge  \frac{2}{3}|V(G)| $.
As before, let $n_0$ be the number  of vertices with $|f(v)|= \emptyset$.
Let us write $n_0 =  n_{01} + n_{02} $,
where  $n_{01}$ is the number of vertices  with exactly one neighbor in $V_2$
and $n_{02}$ is the number of vertices in $V_0$ with more than one neighbor in $V_2$.
Obviously, each vertex in $V_0$ must have at least one neighbor with more than one color.
Observe that 
$$ 2n_{01} \leq  3n_1 ~~~{   \textrm{and }} \quad 2n_{02} + n_{01} \leq  3n_2 \, ,
$$
which implies 
\begin{equation}
 \frac{4}{3} n_{02} +  n_{01} \leq 2n_2 + \frac{1}{2} n_1
\end{equation}
Thus we have
\begin{eqnarray}
 w(f) &=& 2 n_2 + n_1 = 
\frac{2}{3}( 3 n_2 + \frac{3}{2} n_1)  =\\
&=& 
\frac{2}{3}(  n_2 + n_1 + n_0 +  2n_2 + \frac{1}{2} n_1 - n_0)  =\\
&=& 
\frac{2}{3}n +  \frac{2}{3}( 2n_2 + \frac{1}{2} n_1 - n_{01}  - n_{02} )  \geq\\
&\geq& 
\frac{2}{3}n +  \frac{2}{3}( \frac{4}{3} n_{02} +  n_{01} - n_{01}  - n_{02} )  \geq\\
&\geq& 
\frac{2}{3}n
\end{eqnarray}
as needed.
\qed
\end{proof}

\medskip
\noindent{\bf Remark.}   
Proposition \ref{SplosnaMeja4}  is a special case of Theorem  \ref{LBKuzman} (Theorem 1.2 in  \cite{Kuzman}). 
The proof based on Lemma \ref{SplosnaMeja4lema} was provided for completeness of presentation. 
It also may be useful to get some more insight in the structure of the functions   with $ w(f)$ equal to $\gamma_{r4} $ in general case. 

\medskip\noindent
The lower bound for 5RD can be proved similarly.

\begin{lem}  \label{SplosnaMeja5lema} 
 Let $G$ be a 3-regular graph on $n$ vertices and $f$ be a $5RD$ function. Then  $  w(f) <   |V(G)| $  implies  $n_4=n_5=0$. 
\end{lem}   

\begin{proof} 
Assume that  the lemma  does not hold. 
Let $G$ be a minimal counterexample, i.e.  $  w(f) <  \frac{3}{4}|V(G)| $ for some5RDF $f$   with $|V(G)| =n$ and   $n_4\neq 0$ or $n_5 \neq 0$. 
Suppose  there is a vertex $v \in V(G)$ with $w(f(v)) = w_v \in\{4,5\}$
and  construct   $G'$ as   in the proof of Lemma  \ref{SplosnaMeja4lema}. 
Observe that  $w(f(G')) = w(f(G)) - w_v$, and $n' = |V(G')| =n - 4$, so 
$$  w(f(G'))   = w(f(G)) - w_v <  n - w_v =   (n' + 4) - w_v  =   n' +  4 - w_v   <    n'  \, ,
$$
 contradicting the minimality of $G$.
\qed
\end{proof}  

\medskip\noindent
Lemma \ref{SplosnaMeja5lema}, in analogy to Lemma   \ref{SplosnaMeja4lema}, gives information on minimal weight 5RDF's
and allows a direct proof of lower bound  for $ \gamma_{r5} $, which is another special case of Theorem  \ref{LBKuzman} (Theorem 1.2 in  \cite{Kuzman}). 

\begin{proposition}  \label{SplosnaMeja5} 
 Let $G$ be a 3-regular graph on $n$ vertices and $f$ be a $5RD$ function. Then $  w(f) \geq  \frac{5}{6}|V(G)|    \, .$ 
\end{proposition}   

\begin{proof} Assume that $f$ is a 5RDF  of $G$ with  $n_4=n_5=0$.  
There are $n_0$ vertices with no color assigned, and each must have five colors in the neighborhood. 
A vertex $v \in V_1$ satisfies demand one of three neighbors, provided they are all in $V_0$. 
Similarly, a vertex  $v \in V_2$ satisfies demand  two  of three neighbors, provided they are all in $V_0$. 
So we have 
\begin{equation}
  5n_0 \le   3 n_1 + 6 n_2 +  9 n_3 \, .
 \label{eqn0xx}
\end{equation}
For the  weight  of  $f$  we have 
\begin{eqnarray}
w(f) &=& 3n_3 + 2n_2 + n_1 = \frac{5}{6}\left\{   \frac{18}{5}n_3 + \frac{12}{5}n_2 + \frac{6}{5}n_1     \right\} \nonumber \\
       &=& \frac{5}{6}\left\{ n_0 +  n_1 + n_2 +  n_3 + ( \frac{13}{5}n_3 +\frac{7}{5}n_2 + \frac{1}{5}n_1    -  n_0 ) \right\}   \, .
 \label{eqn04xx} 
\end{eqnarray}
From (\ref{eqn0xx})   we   have $  \frac{13}{5}n_3 +\frac{7}{5}n_2 + \frac{1}{5}n_1   > \frac{9}{5}n_3 +\frac{6}{5}n_2 + \frac{3}{5}n_1      \ge  n_0  $
implying 
$  f(w)   \ge   \frac{5}{6}   n   \, ,$
as claimed.
\qed
\end{proof}  

\medskip\noindent
We continue with a more precise analysis of the 4RDF for graphs that attain the lower bound.

\begin{lem} \label{DvodelenGraf4}
Let $G$ be a cubic graph   with  $ \gamma_{r4}(G) =  \frac{2}{3}n$  and let $f$ be a 4RD function of minimal weight, $|f(G)| =  \frac{2}{3}n$. 
Then, exactly one half  of the vertices   are colored.  
Futhermore, the number of vertices that are assigned two colors is 
$c_2  = \frac{n}{6}$ and   the number of vertices that are assigned one  color is  $c_1  = \frac{n}{3}$.
\end{lem}

\begin{proof}  
Let $c_0$ be the number of vertices that are not colored by $f$. 
From $ \gamma_{r4}(G) =  \frac{2}{3}n$   we know that the average weight of  vertices is $\frac{2}{3}$.
We are going to use a standard discharging argument.
Start with a 4DRF $f$, and discharge the vertices as follows. 
Vertices  of weight 2 give   $\frac{4}{9}$  to the three neighbors   and   vertices with one color give 
 $\frac{1}{9}$  to the three neighbors. 
As $f$ is a 4RDF, every uncolored vertex receives at least $\frac{2}{3}$, hence  the total weight transferred is at least $\frac{2}{3} c_0$.
On the other hand, the colored vertices are left with weight   $\frac{2}{3}$ or more, if they receive any charge from a neighbor that is colored by $f$.
As the total weight is $\frac{2}{3}n$, we conclude that all the weights after discharging must be $\frac{2}{3}$.
The colored vertices in total transfer weight $3\frac{4}{9} c_2 +  3  \frac{1}{9}  c_1$ so we have
\begin{equation} 
\frac{2}{3} c_0 = 3\frac{4}{9} c_2 +  3  \frac{1}{9}  c_1 \label{eqn1}
\end{equation}
Furthermore,  we have, by definition 
\begin{equation} 
c_0 + c_1  + c_2 = n  \, , \label{eqn2}
\end{equation} 
and, by assumption 
\begin{equation} 
 c_1  + 2c_2 =\frac{2}{3}n\,. \label{eqn3}
\end{equation}  
The unique solution of the system (\ref{eqn1}), (\ref{eqn2}), and (\ref{eqn3}) is  $c_0  = \frac{n}{2}$,  $c_1  = \frac{n}{3}$, and $c_2  = \frac{n}{6}$, as claimed. 
\qed
\end{proof}  

\medskip\noindent
The proof of  Lemma \ref{DvodelenGraf4} also implies that $G$ must be bipartite, each uncolored vertex has exactly one neighbor colored with two colors and two neighbors 
colored with one color.
This means that the sets of colors of neighbors must form a partition of the color set.
We write these observations formally for later reference.

\begin{lem} \label{BipartiteAndPartition}
Let $G$ be a cubic graph   with  $ \gamma_{r4}(G) =  \frac{2}{3}n$  and let $f$ be a 4RD function of minimal weight, $|f(G)| =  \frac{2}{3}n$. 
Then\\
(1)  $G$ is bipartite, all colored vertices are in one, and all uncolored vertices are in the second set of vertices;\\
(2)  the color sets assigned to neighbors of an uncolored vertex form a partition of the set of colors.
\end{lem}

\medskip\noindent
Similar reasoning for 5RDF is seemingly not sufficient for an analogous result. We  can only prove 

\begin{lem} \label{DvodelenGraf5}
Let $G$ be a cubic graph and with  $ \gamma_{r5}(G) =  \frac{5}{6}n$   
and let $f$ be a 5RD function of minimal weight $|f(G)| =  \frac{5}{6}n$. 
Then, exactly one half  of the vertices   are colored and G is a bipartite graph, where all colored vertices are in one, and all uncolored vertices are in the second set of vertices. 
Furthermore, the neighbors of an  uncolored vertex are colored with color sets that form a partition of the  color set.
\end{lem}

\begin{proof}
Let $f$ be a 5RD function of minimal weight $|f(G)| =  \frac{5}{6}n$ and define charges on vertices with $w(v) = |f(v)|$.
Let the discharging rules be:
vertex with weight 3  gives $ \frac{13}{18}$ to each of the three neighbors;
vertex with weight 2 gives $ \frac{7}{18}$ to each of the three neighbors; and 
vertex with weight 1 gives $ \frac{1}{18}$ to each of the three neighbors.
Observe that the charges of colored vertices, provided they form an independent set, is $ \frac{5}{6}$.
As $f$ is a 5RDF, each uncolored vertex receives (at least)  in total $ \frac{5}{6}$.
Thus, after discharging, all vertices have charge at least  $ \frac{5}{6}$, and since the average charge is     $ \frac{5}{6}$, 
the charges must be exactly   $ \frac{5}{6}$.
It follows that the set of colored vertices is independent.
Furthermore, each uncolored vertex must have three colored neighbors, because
 $ \frac{5}{6}$ can  either  be get as $ \frac{13}{18} + 2 \frac{1}{18} = \frac{15}{18}$, or  $2  \frac{7}{18} +  \frac{1}{18} = \frac{15}{18}$.
(The sum of discharged weights from vertices of weight 2 and 3 is  $ \frac{13}{18} +  \frac{7}{18} >  \frac{5}{6}$.)
Hence each uncolored vertex must have three colored neighbors.
Obviously,  the neighbors of an  uncolored vertex are colored with color sets that form a partition of the  color set.
This concludes the proof.
%
%....................... {\color{red} SISTEM ENAČB  ZA 4 NEZNANKE ................\\
%
%ALI LP....  ...........  PREVERI !!!!  }
\qed
\end{proof}

%%%%%%%%%%%%%%%%%%%%%%%%%%%%%%%%%%%%%%%%%%%%%%%%%%%%%%%%%%
\subsection{General upper bound}
 
Recall that  Proposition  \ref{TheoremGeneral}  holds for RDR graphs.  
In order to use the 3-rainbow domination upper  bounds for  establishing upper  bounds for  4 and 5-rainbow domination, 
we need a slightly weaker statement that holds for all graphs.

\begin{Theorem}  \label{TheoremGeneralNOV}  
Let   $f$  be an $\ell$-RDF  of $G$ with  $w(f) = \gamma_{r\ell}(G)$.
Define $U_i = \{ v ~|~  i \in f(v)  {\rm ~for~} \ell  {\rm ~RDF~} f \}$  so that
$\gamma_{r\ell}(G)  =  |U_1| + |U_2| + \dots  + |U_\ell|$.
It holds
\begin{equation}
  \gamma_{r(\ell+1)}(G)  \le  \gamma_{r\ell}(G)  + \min_i \{ |U_i| \} 
 \le  \gamma_{r\ell}(G)  + \left\lfloor \frac{ \gamma_{r\ell}(G)}{\ell} \right\rfloor   =   \left\lfloor \frac{\ell+1 }{\ell}  \gamma_{r\ell}(G) \right\rfloor 
\label{neenakostTheoremGeneralNOV}
\end{equation}
\end{Theorem} 
 
\begin{proof}
Assume $f$  is an  $\ell$-RDF  of $G$ with  $w(f) = \gamma_{r\ell}(G)$.
We will define a function $f'$ by adding the color $\ell +1$ to some vertices of $G$.

Let  $i \in \{1,2,\dots,\ell\}$ be an arbitrary color.
For any  vertex  $v$    of $G$ with $f(v) = \emptyset$,   $v$ must have a neighbor, say $u$ such that $i\in f(u)$. 
Set $f'(u) =  f(u) \cup \{ \ell +1\}$. % and observe that  $u$ is $\ell +1$ rainbow dominated by $f'$. 
Apply the rule to all vertices  $v$    of $G$ with $f(v) = \emptyset$, and 
define $f'(x) = f(x)$ for all vertices $u$ that were not selected as neighbors of some $v$ with  $f(v) = \emptyset$.
Observe that all vertices are   $\ell +1$ rainbow dominated by $f'$. 
Clearly,  $|U_{\ell +1}| = |\{ v ~|~   (\ell+1) \in f(v)  {\rm ~for~} \ell  {\rm ~RDF~} f \}| \le |U_i|$.
As any $i$ can be chosen, the inequality (\ref{neenakostTheoremGeneralNOV}) follows.
\qed
\end{proof}

\medskip 
\noindent{\bf Remark.}    
When $G$ is $d$-regular graph with   $\gamma_{r\ell}(G)  =\frac{\ell |V(G)|}{2d}$ then 
$  \gamma_{r(\ell+1)}(G) =\frac{(\ell +1)|V(G)|}{2d}$,  
c.f. Proposition \ref{TheoremGeneral} (Proposition 3.4 in \cite{Kuzman}).   
Also note that the statement of  Theorem  \ref{TheoremGeneralNOV}  
may be useful only when   $d \le  \ell <  2d$ as we know that  for $\ell \ge 2d$
 $  \gamma_{r\ell}(G) = |V(G)|$.

%%%%%%%%%%%%%%%%%%%%%%%%%%%%%%%%%%%%%%%%%%%%%%%%%%%%%%%%%%%%%%%%%%%%%%%%%%%%%%%%%%%%%%
\section{Rainbow domination of Petersen graphs}

\subsection{More definitions and previous work}
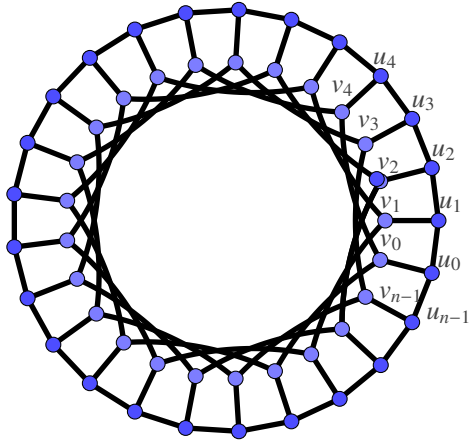
\begin{figure}[h]
	\centering 
	
	\definecolor{uuuuuu}{rgb}{0.26666666666666666,0.26666666666666666,0.26666666666666666}
	\definecolor{oznaka}{rgb}{0.26666666666666666,0.26666666666666666,0.26666666666666666} 
	\definecolor{ududff}{rgb}{0.30196078431372547,0.30196078431372547,1}
	\definecolor{xdxdff}{rgb}{0.49019607843137253,0.49019607843137253,1}
	\begin{tikzpicture}   [scale=0.7,style=thin]
		\def\vr{4pt}
		\draw [line width=2pt] (3,0)-- (4,0);
		\draw [line width=2pt] (2.9057494833858932,0.7460696614945643)-- (3.8743326445145243,0.9947595486594192);
		\draw [line width=2pt] (2.6289200401315904,1.4452610223051456)-- (3.505226720175454,1.927014696406861);
		\draw [line width=2pt] (2.186905882264234,2.0536413177860657)-- (2.9158745096856458,2.7381884237147545);
		\draw [line width=2pt] (2.9057494833858932,-0.7460696614945643)-- (3.8743326445145243,-0.9947595486594192);
		\draw [line width=2pt] (0.9270509831248424,2.8531695488854605)-- (1.2360679774997898,3.804226065180614);
		\draw [line width=2pt] (0.1883715585879402,2.9940801852848145)-- (0.25116207811725355,3.992106913713086);
		\draw [line width=2pt] (-0.5621439437571737,2.9468617521860656)-- (-0.7495252583428986,3.9291490029147544);
		\draw [line width=2pt] (-1.2773378746952178,2.714481157398058)-- (-1.7031171662602904,3.6193082098640774);
		\draw [line width=2pt] (1.6074803849369899,2.532983776506045)-- (2.1433071799159866,3.37731170200806);
		\draw [line width=2pt] (-2.427050983124842,1.7633557568774194)-- (-3.2360679774997894,2.3511410091698925);
		\draw [line width=2pt] (-2.7893294576647536,1.104373658054034)-- (-3.719105943553005,1.472498210738712);
		\draw [line width=2pt] (-2.9763441039434326,0.3759997006929129)-- (-3.968458805257911,0.501332934257217);
		\draw [line width=2pt] (-2.9763441039434326,-0.3759997006929128)-- (-3.96845880525791,-0.5013329342572168);
		\draw [line width=2pt] (-1.9122719692460688,2.3115397283273675)-- (-2.549695958994758,3.082052971103157);
		\draw [line width=2pt] (-2.4270509831248424,-1.7633557568774187)-- (-3.23606797749979,-2.351141009169892);
		\draw [line width=2pt] (-1.912271969246069,-2.3115397283273666)-- (-2.549695958994759,-3.082052971103156);
		\draw [line width=2pt] (-1.277337874695218,-2.714481157398058)-- (-1.7031171662602906,-3.6193082098640774);
		\draw [line width=2pt] (-0.5621439437571738,-2.946861752186065)-- (-0.7495252583428984,-3.929149002914753);
		\draw [line width=2pt] (-2.7893294576647536,-1.1043736580540333)-- (-3.7191059435530054,-1.4724982107387108);
		\draw [line width=2pt] (0.9270509831248416,-2.8531695488854605)-- (1.2360679774997891,-3.8042260651806146);
		\draw [line width=2pt] (1.6074803849369885,-2.532983776506045)-- (2.1433071799159857,-3.3773117020080603);
		\draw [line width=2pt] (2.1869058822642335,-2.0536413177860657)-- (2.9158745096856453,-2.7381884237147545);
		\draw [line width=2pt] (2.6289200401315895,-1.4452610223051456)-- (3.5052267201754526,-1.9270146964068606);
		\draw [line width=2pt] (0.18837155858793975,-2.9940801852848136)-- (0.25116207811725233,-3.9921069137130853);
		\draw [line width=2pt] (3.8743326445145243,0.9947595486594192)-- (4,0);
		\draw [line width=2pt] (3.8743326445145243,0.9947595486594192)-- (3.505226720175454,1.927014696406861);
		\draw [line width=2pt] (3.505226720175454,1.927014696406861)-- (2.9158745096856458,2.7381884237147545);
		\draw [line width=2pt] (2.9158745096856458,2.7381884237147545)-- (2.1433071799159866,3.37731170200806);
		\draw [line width=2pt] (2.1433071799159866,3.37731170200806)-- (1.2360679774997898,3.804226065180614);
		\draw [line width=2pt] (1.2360679774997898,3.804226065180614)-- (0.25116207811725355,3.992106913713086);
		\draw [line width=2pt] (0.25116207811725355,3.992106913713086)-- (-0.7495252583428986,3.9291490029147544);
		\draw [line width=2pt] (-0.7495252583428986,3.9291490029147544)-- (-1.7031171662602904,3.6193082098640774);
		\draw [line width=2pt] (-1.7031171662602904,3.6193082098640774)-- (-2.549695958994758,3.082052971103157);
		\draw [line width=2pt] (-2.549695958994758,3.082052971103157)-- (-3.2360679774997894,2.3511410091698925);
		\draw [line width=2pt] (-3.2360679774997894,2.3511410091698925)-- (-3.719105943553005,1.472498210738712);
		\draw [line width=2pt] (-3.719105943553005,1.472498210738712)-- (-3.968458805257911,0.501332934257217);
		\draw [line width=2pt] (-3.968458805257911,0.501332934257217)-- (-3.96845880525791,-0.5013329342572168);
		\draw [line width=2pt] (-3.96845880525791,-0.5013329342572168)-- (-3.7191059435530054,-1.4724982107387108);
		\draw [line width=2pt] (-3.7191059435530054,-1.4724982107387108)-- (-3.23606797749979,-2.351141009169892);
		\draw [line width=2pt] (-3.23606797749979,-2.351141009169892)-- (-2.549695958994759,-3.082052971103156);
		\draw [line width=2pt] (-2.549695958994759,-3.082052971103156)-- (-1.7031171662602906,-3.6193082098640774);
		\draw [line width=2pt] (-1.7031171662602906,-3.6193082098640774)-- (-0.7495252583428984,-3.929149002914753);
		\draw [line width=2pt] (-0.7495252583428984,-3.929149002914753)-- (0.25116207811725233,-3.9921069137130853);
		\draw [line width=2pt] (0.25116207811725233,-3.9921069137130853)-- (1.2360679774997891,-3.8042260651806146);
		\draw [line width=2pt] (1.2360679774997891,-3.8042260651806146)-- (2.1433071799159857,-3.3773117020080603);
		\draw [line width=2pt] (2.1433071799159857,-3.3773117020080603)-- (2.9158745096856453,-2.7381884237147545);
		\draw [line width=2pt] (2.9158745096856453,-2.7381884237147545)-- (3.5052267201754526,-1.9270146964068606);
		\draw [line width=2pt] (3.5052267201754526,-1.9270146964068606)-- (3.8743326445145243,-0.9947595486594192);
		\draw [line width=2pt] (3.8743326445145243,-0.9947595486594192)-- (4,0);
		\draw [shift={(5.8774236232615715,5.519268155410936)},line width=2pt]  plot[domain=3.635614407983292:4.1555353729193945,variable=\t]({1*5.622656939199773*cos(\t r)+0*5.622656939199773*sin(\t r)},{0*5.622656939199773*cos(\t r)+1*5.622656939199773*sin(\t r)});
		\draw [shift={(4.320187377509104,6.80752601477355)},line width=2pt]  plot[domain=3.8869418202704793:4.406862785206575,variable=\t]({1*5.622656939199847*cos(\t r)+0*5.622656939199847*sin(\t r)},{0*5.622656939199847*cos(\t r)+1*5.622656939199847*sin(\t r)});
		\draw [shift={(2.4914978702899373,7.668041978298544)},line width=2pt]  plot[domain=4.138269232557663:4.658190197493757,variable=\t]({1*5.6226569391998655*cos(\t r)+0*5.6226569391998655*sin(\t r)},{0*5.6226569391998655*cos(\t r)+1*5.6226569391998655*sin(\t r)});
		\draw [shift={(0.5062583887922505,8.046746663241425)},line width=2pt]  plot[domain=4.389596644844851:4.909517609780936,variable=\t]({1*5.6226569391999615*cos(\t r)+0*5.6226569391999615*sin(\t r)},{0*5.6226569391999615*cos(\t r)+1*5.6226569391999615*sin(\t r)});
		\draw [shift={(-1.5107911691614124,7.919844661468734)},line width=2pt]  plot[domain=4.640924057132039:5.1608450220681155,variable=\t]({1*5.622656939200044*cos(\t r)+0*5.622656939200044*sin(\t r)},{0*5.622656939200044*cos(\t r)+1*5.622656939200044*sin(\t r)});
		\draw [shift={(-3.432912161655407,7.295309692464662)},line width=2pt]  plot[domain=4.892251469419222:5.4121724343553,variable=\t]({1*5.6226569392000165*cos(\t r)+0*5.6226569392000165*sin(\t r)},{0*5.6226569392000165*cos(\t r)+1*5.6226569392000165*sin(\t r)});
		\draw [shift={(-5.139330657664762,6.212383585210749)},line width=2pt]  plot[domain=5.1435788817064:5.663499846642488,variable=\t]({1*5.6226569391999215*cos(\t r)+0*5.6226569391999215*sin(\t r)},{0*5.6226569391999215*cos(\t r)+1*5.6226569391999215*sin(\t r)});
		\draw [shift={(-6.522826107317136,4.739110569749559)},line width=2pt]  plot[domain=5.394906293993587:5.914827258929669,variable=\t]({1*5.62265693919999*cos(\t r)+0*5.62265693919999*sin(\t r)},{0*5.62265693919999*cos(\t r)+1*5.62265693919999*sin(\t r)});
		\draw [shift={(-7.496468403370389,2.9680618079614702)},line width=2pt]  plot[domain=5.64623370628077:6.166154671216852,variable=\t]({1*5.622656939199996*cos(\t r)+0*5.622656939199996*sin(\t r)},{0*5.622656939199996*cos(\t r)+1*5.622656939199996*sin(\t r)});
		\draw [shift={(-7.999080019557544,1.0105188070113869)},line width=2pt]  plot[domain=-0.3856241886116374:0.13429677632445355,variable=\t]({1*5.622656939199896*cos(\t r)+0*5.622656939199896*sin(\t r)},{0*5.622656939199896*cos(\t r)+1*5.622656939199896*sin(\t r)});
		\draw [shift={(-7.999080019557695,-1.010518807011407)},line width=2pt]  plot[domain=-0.13429677632444648:0.38562418861163006,variable=\t]({1*5.622656939200043*cos(\t r)+0*5.622656939200043*sin(\t r)},{0*5.622656939200043*cos(\t r)+1*5.622656939200043*sin(\t r)});
		\draw [shift={(-7.49646840337044,-2.9680618079614915)},line width=2pt]  plot[domain=0.11703063596273712:0.6369516008988132,variable=\t]({1*5.622656939200049*cos(\t r)+0*5.622656939200049*sin(\t r)},{0*5.622656939200049*cos(\t r)+1*5.622656939200049*sin(\t r)});
		\draw [shift={(-6.522826107317128,-4.7391105697495535)},line width=2pt]  plot[domain=0.3683580482499173:0.8882790131859998,variable=\t]({1*5.622656939199981*cos(\t r)+0*5.622656939199981*sin(\t r)},{0*5.622656939199981*cos(\t r)+1*5.622656939199981*sin(\t r)});
		\draw [shift={(-5.139330657664807,-6.212383585210805)},line width=2pt]  plot[domain=0.6196854605371014:1.1396064254731828,variable=\t]({1*5.622656939199993*cos(\t r)+0*5.622656939199993*sin(\t r)},{0*5.622656939199993*cos(\t r)+1*5.622656939199993*sin(\t r)});
		\draw [shift={(-3.4329121616554126,-7.295309692464681)},line width=2pt]  plot[domain=0.8710128728242871:1.3909338377603642,variable=\t]({1*5.622656939200035*cos(\t r)+0*5.622656939200035*sin(\t r)},{0*5.622656939200035*cos(\t r)+1*5.622656939200035*sin(\t r)});
		\draw [shift={(-1.5107911691614127,-7.919844661468742)},line width=2pt]  plot[domain=1.1223402851114712:1.6422612500475469,variable=\t]({1*5.622656939200053*cos(\t r)+0*5.622656939200053*sin(\t r)},{0*5.622656939200053*cos(\t r)+1*5.622656939200053*sin(\t r)});
		\draw [shift={(0.5062583887922607,-8.046746663241567)},line width=2pt]  plot[domain=1.373667697398657:1.893588662334728,variable=\t]({1*5.6226569392001*cos(\t r)+0*5.6226569392001*sin(\t r)},{0*5.6226569392001*cos(\t r)+1*5.6226569392001*sin(\t r)});
		\draw [shift={(2.4914978702899884,-7.668041978298696)},line width=2pt]  plot[domain=1.624995109685837:2.144916074621915,variable=\t]({1*5.622656939200022*cos(\t r)+0*5.622656939200022*sin(\t r)},{0*5.622656939200022*cos(\t r)+1*5.622656939200022*sin(\t r)});
		\draw [shift={(4.320187377509192,-6.807526014773685)},line width=2pt]  plot[domain=1.8763225219730195:2.3962434869090994,variable=\t]({1*5.622656939200005*cos(\t r)+0*5.622656939200005*sin(\t r)},{0*5.622656939200005*cos(\t r)+1*5.622656939200005*sin(\t r)});
		\draw [shift={(5.8774236232617145,-5.519268155411066)},line width=2pt]  plot[domain=2.1276499342602015:2.6475708991962854,variable=\t]({1*5.622656939199962*cos(\t r)+0*5.622656939199962*sin(\t r)},{0*5.622656939199962*cos(\t r)+1*5.622656939199962*sin(\t r)});
		\draw [shift={(7.0653597271126545,-3.8842143793956194)},line width=2pt]  plot[domain=2.3789773465473836:2.8988983114834697,variable=\t]({1*5.622656939199937*cos(\t r)+0*5.622656939199937*sin(\t r)},{0*5.622656939199937*cos(\t r)+1*5.622656939199937*sin(\t r)});
		\draw [shift={(7.809353294733635,-2.005101128781521)},line width=2pt]  plot[domain=2.630304758834564:3.150225723770656,variable=\t]({1*5.622656939199871*cos(\t r)+0*5.622656939199871*sin(\t r)},{0*5.622656939199871*cos(\t r)+1*5.622656939199871*sin(\t r)});
		\draw [shift={(8.062656474054178,0)},line width=2pt]  plot[domain=2.8816321711217467:3.4015531360578404,variable=\t]({1*5.622656939199856*cos(\t r)+0*5.622656939199856*sin(\t r)},{0*5.622656939199856*cos(\t r)+1*5.622656939199856*sin(\t r)});
		\draw [shift={(7.809353294733636,2.0051011287815284)},line width=2pt]  plot[domain=3.132959583408931:3.6528805483450233,variable=\t]({1*5.622656939199874*cos(\t r)+0*5.622656939199874*sin(\t r)},{0*5.622656939199874*cos(\t r)+1*5.622656939199874*sin(\t r)});
		\draw [shift={(7.809353294733635,-2.005101128781521)},line width=2pt]  plot[domain=2.630304758834564:3.150225723770656,variable=\t]({1*5.622656939199871*cos(\t r)+0*5.622656939199871*sin(\t r)},{0*5.622656939199871*cos(\t r)+1*5.622656939199871*sin(\t r)});
		\draw [shift={(7.065359727112454,3.884214379395512)},line width=2pt]  plot[domain=3.384286995696106:3.904207960632214,variable=\t]({1*5.6226569391997145*cos(\t r)+0*5.6226569391997145*sin(\t r)},{0*5.6226569391997145*cos(\t r)+1*5.6226569391997145*sin(\t r)});
		%\begin{scriptsize}
			\draw [fill=xdxdff] (3,0) circle (\vr);
			\draw[color=oznaka] (3.1112343865503296,0.3034062300017742) node {$v_1$};
			\draw [fill=ududff] (4,0) circle (\vr);
			\draw[color=oznaka] (4.115687460821578,0.3034062300017742) node {~~$u_1$};
			%\draw [fill=uuuuuu] (0,0) circle (2pt);
			%\draw[color=uuuuuu] (0.11182590087923941,0.2755047557164618) node {$C$};
			\draw [fill=xdxdff] (2.9057494833858932,0.7460696614945643) circle (\vr);
			\draw[color=oznaka] (3.0693821751223602,1.0427952985625528) node {$v_2$};
			\draw [fill=ududff] (3.8743326445145243,0.9947595486594192) circle (\vr);
			\draw[color=oznaka] (4.045933775108296,1.2939085671303643) node {~$u_2$};
			\draw [fill=xdxdff] (2.6289200401315904,1.4452610223051456) circle (\vr);
			\draw[color=oznaka] (2.7,1.85) node {$v_3$};
			\draw [fill=ududff] (3.505226720175454,1.927014696406861) circle (\vr);
			\draw[color=oznaka] (3.7250668208272035,2.2286079556883296) node {$u_3$};
			\draw [fill=xdxdff] (2.186905882264234,2.0536413177860657) circle (\vr);
			\draw[color=oznaka] (2.2,2.5) node {$v_4$};
			\draw [fill=ududff] (2.9158745096856458,2.7381884237147545) circle (\vr);
			\draw[color=oznaka] (3,3.0377507099623893) node {$u_4$};
			\draw [fill=xdxdff] (2.9057494833858932,-0.7460696614945643) circle (\vr);
			\draw[color=oznaka] (3.111234386550329,-0.4080813642736919) node {$v_0$};
			\draw [fill=ududff] (3.8743326445145243,-0.9947595486594192) circle (\vr);
			\draw[color=oznaka] (4.087785986536265,-0.8) node {~~$u_0$};
			\draw [fill=xdxdff] (0.9270509831248424,2.8531695488854605) circle (\vr);
			%\draw[color=xdxdff] (1.1302297122931444,3.1912088185316074) node {$A'_{2}$};
			\draw [fill=ududff] (1.2360679774997898,3.804226065180614) circle (\vr);
			%\draw[color=ududff] (1.4510966665742377,4.139858944232229) node {$B'_{2}$};
			\draw [fill=xdxdff] (0.1883715585879402,2.9940801852848145) circle (\vr);
			%\draw[color=xdxdff] (0.446643592302989,3.3307161899581694) node {$A''_{1}$};
			\draw [fill=ududff] (0.25116207811725355,3.992106913713086) circle (\vr);
			%\draw[color=ududff] (0.5163972780162701,4.335169264229416) node {$B''_{1}$};
			\draw [fill=xdxdff] (-0.5621439437571737,2.9468617521860656) circle (\vr);
			%\draw[color=xdxdff] (-0.2369425276871664,3.288863978530201) node {$A'''_{1}$};
			\draw [fill=ududff] (-0.7495252583428986,3.9291490029147544) circle (\vr);
			%\draw[color=ududff] (-0.4322528476843537,4.2654155785161345) node {$B'''_{1}$};
			\draw [fill=xdxdff] (-1.2773378746952178,2.714481157398058) circle (\vr);
			%\draw[color=xdxdff] (-1.171641916245134,3.009849235677077) node {$D$};
			\draw [fill=ududff] (-1.7031171662602904,3.6193082098640774) circle (\vr);
			%\draw[color=ududff] (-1.590164030524821,3.9166471499497297) node {$E$};
			\draw [fill=xdxdff] (1.6074803849369899,2.532983776506045) circle (\vr);
			%\draw[color=xdxdff] (1.716160672284706,2.8284896528225465) node {$F$};
			\draw [fill=ududff] (2.1433071799159866,3.37731170200806) circle (\vr);
			%\draw[color=ududff] (2.260239420848299,3.6794846185245746) node {$G$};
			\draw [fill=xdxdff] (-2.427050983124842,1.7633557568774194) circle (\vr);
			%\draw[color=xdxdff] (-2.3156023619429447,2.0611991099764553) node {$H$};
			\draw [fill=ududff] (-3.2360679774997894,2.3511410091698925) circle (\vr);
			%\draw[color=ududff] (-3.124745116217006,2.6471300699680156) node {$I$};
			\draw [fill=xdxdff] (-2.7893294576647536,1.104373658054034) circle (\vr);
			%\draw[color=xdxdff] (-2.678321527652007,1.4055144642716138) node {$J$};
			\draw [fill=ududff] (-3.719105943553005,1.472498210738712) circle (\vr);
			%\draw[color=ududff] (-3.6130209162099742,1.7682336299806751) node {$K$};
			\draw [fill=xdxdff] (-2.9763441039434326,0.3759997006929129) circle (\vr);
			%\draw[color=xdxdff] (-2.859681110506538,0.6800761328534917) node {$L$};
			\draw [fill=ududff] (-3.968458805257911,0.501332934257217) circle (\vr);
			%\draw[color=ududff] (-3.85018344763513,0.8056327671373973) node {$M$};
			\draw [fill=xdxdff] (-2.9763441039434326,-0.3759997006929128) circle (\vr);
			%\draw[color=xdxdff] (-2.803878161935913,-0.07326367284994317) node {$D'$};
			\draw [fill=ududff] (-3.96845880525791,-0.5013329342572168) circle (\vr);
			%\draw[color=ududff] (-3.7943804990645056,-0.19882030713384896) node {$E'$};
			\draw [fill=xdxdff] (-1.9122719692460688,2.3115397283273675) circle (\vr);
			%\draw[color=xdxdff] (-1.7436221390940394,2.605277858540047) node {$F'$};
			\draw [fill=ududff] (-2.549695958994758,3.082052971103157) circle (\vr);
			%\draw[color=ududff] (-2.3853560476562263,3.3865191385287945) node {$G'$};
			\draw [fill=xdxdff] (-2.4270509831248424,-1.7633557568774187) circle (\vr);
			%\draw[color=xdxdff] (-2.25979941337232,-1.4683373871155632) node {$H'$};
			\draw [fill=ududff] (-3.23606797749979,-2.351141009169892) circle (\vr);
			%\draw[color=ududff] (-3.0689421676463815,-2.054268347107123) node {$I'$};
			\draw [fill=xdxdff] (-1.912271969246069,-2.3115397283273666) circle (\vr);
			%\draw[color=xdxdff] (-1.7436221390940394,-2.0124161356791546) node {$J'$};
			\draw [fill=ududff] (-2.549695958994759,-3.082052971103156) circle (\vr);
			%\draw[color=ududff] (-2.3853560476562263,-2.7797066785252458) node {$K'$};
			\draw [fill=xdxdff] (-1.277337874695218,-2.714481157398058) circle (\vr);
			%\draw[color=xdxdff] (-1.1158389676745089,-2.4169875128161844) node {$L'$};
			\draw [fill=ududff] (-1.7031171662602906,-3.6193082098640774) circle (\vr);
			%\draw[color=ududff] (-1.534361081954196,-3.323785427088837) node {$M'$};
			\draw [fill=xdxdff] (-0.5621439437571738,-2.946861752186065) circle (\vr);
			%\draw[color=xdxdff] (-0.33459768768576004,-2.640199307098684) node {$D''$};
			\draw [fill=ududff] (-0.7495252583428984,-3.929149002914753) circle (\vr);
			%\draw[color=ududff] (-0.5299080076829473,-3.630701644227274) node {$E''$};
			\draw [fill=xdxdff] (-2.7893294576647536,-1.1043736580540333) circle (\vr);
			%\draw[color=xdxdff] (-2.5667156305107572,-0.7987020042680655) node {$F''$};
			\draw [fill=ududff] (-3.7191059435530054,-1.4724982107387108) circle (\vr);
			%\draw[color=ududff] (-3.5014150190687245,-1.175371907119783) node {$G''$};
			\draw [fill=xdxdff] (0.9270509831248416,-2.8531695488854605) circle (\vr);
			%\draw[color=xdxdff] (1.1441804494358006,-2.5564948842427464) node {$H''$};
			\draw [fill=ududff] (1.2360679774997891,-3.8042260651806146) circle (\vr);
			%\draw[color=ududff] (1.465047403716894,-3.505145009943368) node {$I''$};
			\draw [fill=xdxdff] (1.6074803849369885,-2.532983776506045) circle (\vr);
			\draw [fill=ududff] (2.1433071799159857,-3.3773117020080603) circle (\vr);
			\draw [fill=xdxdff] (2.1869058822642335,-2.0536413177860657) circle (\vr);
			\draw[color=oznaka] (3.3,-1.4) node {$v_{n-1}$};
			\draw[color=oznaka] (4.2,-1.8) node {$u_{n-1}$};
			\draw [fill=ududff] (2.9158745096856453,-2.7381884237147545) circle (\vr);
			\draw [fill=xdxdff] (2.6289200401315895,-1.4452610223051456) circle (\vr);
			\draw [fill=ududff] (3.5052267201754526,-1.9270146964068606) circle (\vr);
			\draw [fill=xdxdff] (0.18837155858793975,-2.9940801852848136) circle (\vr);
			\draw [fill=ududff] (0.25116207811725233,-3.9921069137130853) circle (\vr);
			\draw [fill=ududff] (2.841164006743404,0.7868904658100392) circle (\vr);
			
		\end{tikzpicture} 
 \caption{A generalized Petersen graph $P(n,k)$ }\label{FigPetersen}
\end{figure}

For $n\geq 3$ and    $k$,  $1 \leq k \leq  n-1$, 
the generalized Petersen graph $P(n,k)$,   is a graph on $2n$ vertices with
$V(P(n,k))=\{v_i,u_i \mid 0 \leq  i \leq n-1\}$ and edges
%$E(P(n,k))=\{\{u_i,u_{i+1}\}, \{u_i,v_i\}, \{v_i,v_{i+k}\} \mid 0 \leq  i \leq n-1\}$. 
$E(P(n,k))=\{u_iu_{i+1}, u_iv_i, v_iv_{i+k} \mid 0 \leq  i \leq n-1,\}$. 
where   all subscripts are taken modulo $n$.
This standard notation was introduced by Watkins~\cite{Watkins}
(see Figure~\ref{FigPetersen}). 
For convenience, throughout the paper, all subscripts will be taken modulo $n$.
%Clearly, the set of vertices $U= \{u_i \mid 0 \leq  i \leq n-1\}$ 
%induces a cycle that is called the {\em outer cycle}, 
%and when $n=ck$, the set of vertices $V= \{v_i \mid 0 \leq  i \leq n-1\}$ 
%induces $k$ cycles   called the {\em inner cycles}.

Generalized Petersen graphs $P(n,k)$ are 3-regular unless $k=\frac{n}{2}$.
It is well known that  Petersen graphs  $P(n,k)$ and $P(n,n-k)$ are isomorphic, so it is natural to restrict attention  to $P(n,k)$
with $n\geq 3$ and  $k$,  $1 \leq k < \frac{n}{2}$.     

%Petersen graphs are cubic graphs. 
Rainbow domination of $P(ck,k)$ for $r=1,2$ and 3 has been studied recently. 
Here we recall and summarize the main results. 
Note that for  $k=1$  $P(ck,k) = P(k,1) $ is isomorphic to the Cartesian product $C_k \Box K_2$.
To avoid trivialities,  some previous studies have  assumed $k>3$.  

The following results  from \cite{Darja3RD} 
that are generalising the study   \cite{Erves2021b}  of rainbow domination of  $ P(6k,k) $ 
are  of particular importance for present work.
Recall that singleton 3-rainbow domination number   $\tilde\gamma_{r3}$  is the minimum weight over
 3-rainbow domination functions that assign singleton sets (and empty set).
Bounds for 3-rainbow domination of generalized Petersen graphs $P(ck,k)$  that were proven  in \cite{Darja3RD} are   
summarized in Theorems \ref{MainTheoremOLD}, \ref{MainTheorem2OLD}, and \ref{MainTheorem3OLD}.

 \begin{Theorem}  \label{MainTheoremOLD}  \cite{Darja3RD} 
 Let $c \equiv 0\pmod 6$.
 Then,  for   
 3-rainbow domination number $\gamma_{r3}$
 and 
 singleton 3-rainbow domination number 
 $\tilde\gamma_{r3}$
 of generalized Petersen graphs $P(ck,k)$ 
  it holds:  %we~have the following:
 \begin{itemize}
 \item 
  If $ k \equiv 1,5\pmod 6$, then 
  $ \gamma_{r3}(P(ck,k))  =\tilde\gamma_{r3}(P(ck,k))  =ck$;
 \item 
  If $ k \equiv 0\pmod 2$, then $
  ck < \gamma_{r3}(P(ck,k))  \leq \tilde\gamma_{r3}(P(ck,k))  =  c (k+ \frac{1}{2})$;
 \item 
  If $k \equiv 3\pmod 6$,  then  
 $ck   < \gamma_{r3}(P(ck,k))  \leq \tilde\gamma_{r3}(P(ck,k))  \leq c(k+1)$. 
 \end{itemize}
\end{Theorem} 
 
 \begin{Theorem}  \label{MainTheorem2OLD} \cite{Darja3RD} 
 Let $c$ be odd. 
 Then, for 
 3-rainbow domination number $\gamma_{r3}$
 and 
 singleton 3-rainbow domination number 
 $\tilde\gamma_{r3}$
 of generalized Petersen graphs $P(ck,k)$ 
   we~have:
 \begin{itemize}
 \item 
  { 
  If $ k \equiv 1,5\pmod 6$,} then 
  $ ck <
  \gamma_{r3}(P(ck,k)) \leq \tilde\gamma_{r3}(P(ck,k)) = ck + \lceil\frac{k}{2} \rceil$;
 \item 
  { If $ k \equiv 0\pmod 2$,} then $ck <\gamma_{r3}(P(ck,k)) 
       \leq \tilde\gamma_{r3}(P(ck,k)) \leq ck+ 
       \JZ{\lfloor\frac{c}{2}\rfloor}
        + \frac{k}{2}$;
 \item 
  {  
  If $k \equiv 3\pmod 6$,}  then  
 $ck   < \gamma_{r3}(P(ck,k))  \leq \tilde\gamma_{r3}(P(ck,k)) 
 \leq c(k+1) + \lceil \frac{ k-2}{2} \rceil$. 
 \end{itemize}
\end{Theorem} 
 
 \begin{Theorem}  \label{MainTheorem3OLD} \cite{Darja3RD} 
 Let $c$ be even, and $c \not\equiv 0\pmod 6$.
 Then, for 
 3-rainbow domination number $\gamma_{r3}$
 and 
 singleton 3-rainbow domination number 
 $\tilde\gamma_{r3}$
 of generalized Petersen graphs $P(ck,k)$ 
   we~have:
 \begin{itemize}
 \item 
  { 
  If $ k \equiv 1,5\pmod 6$,} then 
    $ ck <  
  \gamma_{r3}(P(ck,k)) \leq \tilde\gamma_{r3}(P(ck,k)) \leq ck + k + 1 $;
 \item 
  { If $ k \equiv 0\pmod 2$,} then $
  ck < \gamma_{r3}(P(ck,k)) \leq \tilde\gamma_{r3}(P(ck,k)) \leq c k+ \frac{c}{2}  + k$;
 \item 
  {  
  If $k \equiv 3~\pmod 6$,}  then  
 $ck   < \gamma_{r3}(P(ck,k))  \leq \tilde\gamma_{r3}(P(ck,k))  \leq ck+  c + k-2$. 
 \end{itemize}
\end{Theorem}   

%%%%%%%%%%%%%%%%%%%%%%%%%%%%%%%%%%%%%%%%%%%%%%%%%%%%%%%%%%%%%%%%%%%%%%%

Theorems \ref{MainTheoremOLD}, \ref{MainTheorem2OLD}, and \ref{MainTheorem3OLD}
imply the following two statements, one covering the cases where exact values are known, and the 
second giving tight lower and  a bit relaxed upper bounds that hold for all other cases.
These bounds will be generalized  to $t=4$ and 5  in Section       \ref{REFTOSECTION}.
%
%First,  let us recall the characterization of graphs $P(ck,k)$   for which  three-rainbow domination is known 
%and equals the general lower bound. This ia part of statement of Theorem 3 in  \cite{Darja3RD} 
\begin{Theorem} \label{DarjaTheorem} %\cite{Darja3RD}  
%Let $k>3$.
% Then if $c\equiv 0~\bmod 6$,   and $k \equiv 1, 5~\bmod 6$,  
%\begin{equation}
%\gamma_{r3}(P(ck, k))   =  ck  \,.
%\end{equation}
$ \gamma_{r3}(P(ck, k))   =  ck$ if and only if  $c\equiv 0~\pmod 6$,   and $k \equiv 1, 5~\pmod 6$.
\end{Theorem} 
%In the same paper, lower and upper bounds are provided in Theorems 4 and 5  in \cite{Darja3RD} . 
%Here we summarize the results in one statement
\begin{Theorem}   \label{MainTheorem3} 
%Let $k>3$. Then 
%\begin{equation}
% ck   \le \gamma_{r3}(P(ck, k))    \le   (c+1)(k+1)  - 3  \,.
%\end{equation}
For $k\ge 1$ and $c\ge 3$ it holds $ck   \le \gamma_{r3}(P(ck, k))    \le   (c+1)(k+1)  - 3  \,.$
\end{Theorem} 

\medskip\noindent
Regarding the 2-rainbow domination, analogous results are summarized in the next two theorems
(for proofs see \cite{Simon2RD}).
%
% \begin{Theorem}  \cite{Simon2RD}    \label{MainTheorem2} 
%For  $c\ge 3$ it holds 
%\begin{equation}
%     \frac{4}{5}ck   \le \gamma_{r2}(P(ck, k))      \le      \frac{4}{5}(c+1)(k+1)  + 1  \,.
%\end{equation}
%\end{Theorem} 
%This is a short and in some cases not best possible bound, as it summarizes a detailed list of upper bounds given in 
%Theorem 4 in \cite{Simon2RD}. Let us only recall the case in which exact value is known.
%\begin{Theorem} %\cite{Simon2RD}  
%Let $k>3$. Then if $c\equiv 0~\bmod 5$,   and $k \equiv 2, 8~\bmod 10$,  
%\begin{equation}
%\gamma_{r2}(P(ck, k))   =     \frac{4}{5}ck  \,.
%\end{equation}
%\end{Theorem} 
\begin{Theorem}  %\cite{Simon2RD} 
   \label{MainTheorem2} 
For $k\ge 1$ and $c\ge 3$ it holds
 $    \frac{4}{5}ck   \le \gamma_{r2}(P(ck, k))      \le      \frac{4}{5}(c+1)(k+1)  + 1  \,.$ 
\end{Theorem} 
\begin{Theorem}  %\cite{Simon2RD} 
$ \gamma_{r2}(P(ck, k))   =  ck$ if and only if   $c\equiv 0~\pmod 5$,   and $k \equiv 2, 8~\pmod {10}$.
\end{Theorem} 

\medskip\noindent
As 1-rainbow domination is just the ordinary domination, we have the following characterization
that holds for  all $P(n,k)$  (compare   \cite{Ebrahimi}, Theorem 1).
\begin{Theorem} \cite{Ebrahimi}  
$\gamma_{r1}(P(n, k))   = \gamma(P(n, k))   =     \frac{1}{2}n$   if and only if 
$n\equiv 0~\pmod 4$, and k odd. 
\end{Theorem} 

\medskip\noindent
In general case, bounds can be extracted from \cite{WangDOMINATION} and summarized as follows.
%\begin{Theorem}\cite{WangDOMINATION}    \label{MainTheorem1} 
%%Let $k>3$. Then if $c\equiv 0~\bmod 5$,   and $k \equiv 2, 8~\bmod 10$,  
%\begin{equation}   
%\frac{1}{2}ck   \le \gamma_{r1}(P(ck, k))      \le      \frac{1}{2}(c+1)k  + 1  \,.
%\end{equation}
%\end{Theorem} 
 \begin{Theorem}\cite{WangDOMINATION}    \label{MainTheorem1}  
For $k\ge 1$ and $c\ge 3$ it holds 
$\frac{1}{2}ck   \le \gamma_{r1}(P(ck, k))      \le      \frac{1}{2}(c+1)k  + 1  \,.$
\end{Theorem} 

%%%%%%%%%%%%%%%%%%%%%%%%%%%%%%%%%%%%%%%%%%%%%%%%%%%%%%%%%%%%%%%%%%%%%%%%%%%%%%%%%%%%%%%%%%%
\subsection{New upper  bounds for $\gamma_{r4}(P(n,k))$ and  $\gamma_{r5}(P(n,k))$.}  

\label{REFTOSECTION}
Here  we apply Theorems \ref{MainTheorem3}  and  \ref{TheoremGeneralNOV}, 
 to obtain  bounds on 4 and 5-rainbow domination numbers of generalized Petersen graphs. 

\begin{Theorem}    \label{MainTheorem4} 
Let  $k\ge 1$ and $c\ge 3$. Then $\gamma_{r4}(P(ck, k))    \le   \frac{4}{3} (c+1)(k+1)  - 4  \,.$
%\begin{equation}
% \frac{4}{3}ck   \le \gamma_{r4}(P(ck, k))    \le   \frac{4}{3} (c+1)(k+1)  - 4  \,.
%\end{equation}
\end{Theorem} 
\begin{proof} 
%For the upper bound, use Theorems \ref{MainTheorem3}  and  \ref{TheoremGeneralNOV}, 
\begin{equation} 
 \gamma_{r4}(P(ck, k))    \le   \left\lfloor \frac{4 }{3 }  \gamma_{r3}(G) \right\rfloor  \le \frac{4}{3} ( (c+1)(k+1)  - 3)    =  \frac{4}{3} (c+1)(k+1)  - 4  \,.
\end{equation}
\qed
\end{proof}

\begin{Theorem}    \label{MainTheorem5} 
Let  $k\ge 1$ and $c\ge 3$. Then $ \gamma_{r5}(P(ck, k))    \le   \frac{5}{3} (c+1)(k+1)  - 5  \,.$
%\begin{equation}
% \frac{5}{3}ck   \le \gamma_{r5}(P(ck, k))    \le   \frac{5}{3} (c+1)(k+1)  - 5  \,.
%\end{equation}
\end{Theorem}  
\begin{proof}
Similarly as in the proof of Theorem   \ref{MainTheorem4},  
 the upper bound follows from  Theorem \ref{MainTheorem4}  using  Theorem  \ref{TheoremGeneralNOV}:
\begin{equation} 
 \gamma_{r5}(P(ck, k))    \le   \left\lfloor \frac{5 }{4 }  \gamma_{r4}(G) \right\rfloor  \le \frac{5}{4} ( \frac{4}{3}(c+1)(k+1)  - 4)    = \frac{5}{3} (c+1)(k+1)  - 5  \,.
\end{equation}
\qed
\end{proof}

Recall that,  for $\ell \ge 6$,  $ \gamma_{r\ell}(P(ck, k))   =   2n$ by Theorem    \ref{LBKuzman}  (Theorem 1.2   in \cite{Kuzman}).
Thus we have tight lower bounds and reasonable upper bounds for  $\gamma_{rt}(P(ck, k)) $,   for all relevant $t$.
 
\subsection{Improved upper bounds}

 The upper bounds can be easily improved by application of   Theorem  \ref{TheoremGeneralNOV} to each of the special cases that are
outlined in 
Theorems \ref{MainTheoremOLD}, \ref{MainTheorem2OLD}, and \ref{MainTheorem3OLD}. 
As the results would most probably not provide exact values, we are not elaborating the details here. 
Also see a comment in concluding section.

%%%%%%%%%%%%%%%%%%%%%%%%%%%%%%%%%%%%%%%%%%%%%%%%%%%%%%%%%%%%%%%%%%%%%%%%%%%%%%%%%%%%%%%%%%%
\subsection{Lower bounds for $\gamma_{r4}(P(n,k))$  and  $\gamma_{r5}(P(n,k))$   - extremal examples.} 

We continue with consideration of generalized Petersen graphs that attain the lower bounds.

\begin{lem} \label{VzorecNaZunanjemCiklu}
Let  $f$ be a rainbow domination function of a generalized Petersen graph $P(n,k)$ with the property that 
each uncolored vertex has its three neighbors colored with three nonempty sets that are a partition of the color set. 
Then, exactly one half  of the vertices on the outer  cycle $\{u_0,u_1,\dots u_{n-1}\}$
are colored. Wlog assume that these are vertices with even indices. 
Then,  the following holds:\\
(1)  $f(u_{i }) =   \emptyset$,  
for~all odd $i$;\\
(2)  $f(u_{i+6}) = f(u_i)$,  $0\leq i \leq  n-1$; %%%%%%%%%%%%%%%%%%%$i\in [0,ck-1]$,
and \\
(3) $f(u_0)$, $f(u_2)$, $f(u_4)$   form a partition of the color set.\\
Consequently, $n\equiv 0\pmod 6$. 
\end{lem}

\medskip\noindent
\begin{proof}  
First observe that the assumed property directly implies  that the graph is bipartite.
Hence we may assume that on the outer cycle exactly one half of the vertices, 
wlog those with even  indices, are colored. 
Then, (1) is trivial.

Another simple but useful observation is that vertices with two consecutive even indices must be colored 
differently, because if $f(u_0) = f(u_2)$, $u_1$ again would not have all three colors in the neighborhood.

Now we   prove   (3). 
Let  the neighbors of $u_1$ be  assigned 
$f(u_0) = A$, $f(u_2) = B$, and $f(v_2)= C$, where nonempty sets $A$, $B$, and $C$ are partition of the color set.
Choose arbitrary colors $a\in A$,  $b\in B$ and $c\in C$.
We claim that  $c\in f(u_4)$.

Assume, for the sake of contradiction, that this  does not hold,  so  that  $c\not \in f(u_4)$.
Then it follows that for the  third neighbor of $u_{3}$ we have $c \in f(v_{3})$.
Also, we know that $c\not\in f(u_{k+1}) $ since $u_{k+1}$ already has one neighbor, $v_{1}$, colored by $c$.
Similarly, $c\not\in f(u_{k+3})  $, as  $u_{3}$ already has one neighbor, $v_{3}$, colored by $c$.
Thus we know that   $c$  is not assigned to   two neighbors of  $ u_{k+2}$. 
Consequently, for the third neighbor of  $u_{k+2}$, we must have $c\in f(v_{2+k })$.

The same reasoning leads to conclusion that $c\in f(v_{2-k })$.
But then, vertex $v_2$ has two neighbors colored by $c$, contradicting the assumption that the color sets of neighbors form a partition of the color set.

As the choice of colors $a$, $b$, and $c$ was arbitrary, it follows that $f(u_4) = C$.
Hence (3) is proven.

Repeating the argument on $u_2$, $u_4$, and  $u_0$, extends the known coloring of the outer cycle to 
 the  local pattern  $A-0-B-0-C-A$.
Straightforward inductive argument completes the proof of (2). 
\qed
\end{proof}  

\medskip\noindent{\bf Remark.} The~lemma says that the vertices on the outer cycle are colored
following the pattern $A-0-B-0-C-0-\dots -A-0-B-0-C-0$,  where  $A$, $B$, and $C$ from a partition of the set of colors
that appears in the neighborhood of an arbitrarily chosen vertex. 

\medskip\noindent 
From the proof of Lemma \ref{VzorecNaZunanjemCiklu} it follows that the length of the outer cycle $n$ must be divisible by 6. 
We continue with observation that the coloring of the outer cycle completely determines the coloring of all vertices, which 
implies neccesary and sufficient conditions for $k$. 
Thus we obtain  a characterization of generalized Petersen graphs with $ \gamma_{r4}(P(n,k)) =  \frac{4}{3}n$ and with $ \gamma_{r5}(P(n,k)) =  \frac{5}{3}n$.
Formally : 

\begin{lem} \label{VzorecNaNotranjemCiklu4}
Assume  $ \gamma_{r4}(P(n,k)) =  \frac{4}{3}n$.  Then,  $n \equiv 0\pmod 6$, and $k \bmod 6$ must be $1$ or $5$.
\end{lem}

\begin{proof} 
Recall that  by Lemma \ref{BipartiteAndPartition},  $ \gamma_{r4}(P(n,k)) =  \frac{4}{3}n$ implies that the corresponding 4RDF 
has the property assumed in  Lemma \ref{VzorecNaZunanjemCiklu}. 
Hence $n \equiv 0\pmod 6$.

Observe that the coloring of the outer cycle uniquely determines the coloring of $P(n,k)$. 
Let $f$ be a 3RD function of minimal weight, $|f(P(n,k))| = n$. 
For example,  by Lemma \ref{VzorecNaZunanjemCiklu},  $f(u_0) = A$, $f(u_2) = B$, and $f(u_4)   = C$
implies $f(u_{6i}) = A$, $f(u_{6i+2}) = B$, and $f(u_{6+4i})   = C$, which in turn 
implies  $f(v_{6i+1}) = C$, $f(v_{6i+3}) = A$, and $f(v_{6i+5}) = C$.
It follows that $v_{6i}$ has a neighbor colored by $A$, and must have two more colored neighbors, hence $k$ must be odd. 
Also, $v_{6i}$ may not be connected to any vertices that already have any color from $A$, hence $k\not= 3 \pmod 6$.
\qed
\end{proof}  

\medskip\noindent
Identical reasoning gives.

\begin{lem} \label{VzorecNaNotranjemCiklu5}
Assume  $ \gamma_{r5}(P(n,k)) =  \frac{5}{3}n$.  Then, $n \equiv 0\pmod 6$, and $k \bmod 6$ must be $1$ or $5$.
\end{lem}

\begin{proof}  Ommited.
\qed
\end{proof}  

\medskip\noindent{\bf Remark.} 
Lemma \ref{VzorecNaZunanjemCiklu} implies that any 4RDF and any 5RDF have  essentially the same structure.
Namely, it is obtained from one of  (the 12 different) 3RDFs by replacing the singleton sets with 
a partition of the color set. In particular, a 5RDF may be constructed by adding two colors to one of the singletons, 
or by enlarging two singleton sets to pairs.

%\medskip\noindent{\bf Remark.} 
Furthermore,  by Proposition \ref{TheoremGeneral},
$ \gamma_{r3}(P(n,k)) =  n$  implies  $ \gamma_{r4}(P(n,k)) =  \frac{4}{3}n$, and   $ \gamma_{r5}(P(n,k)) =  \frac{5}{3}n$. 
So the implications in Lemma(s)  \ref{VzorecNaNotranjemCiklu4} and \ref{VzorecNaNotranjemCiklu5} are equivalences.  
Thus the characterizations of minimal examples among the generalized Petersen graphs 
are equivalent for $r=3$, $4$, and $5$.
This is formally stated as Theorem  \ref{MainTheoremX2novi}  in the first section.

\medskip\noindent
This opens the following questions:

\medskip\noindent{\bf Question 1.}
Are there cubic graphs with  $ \gamma_{r4}(P(n,k)) =  \frac{4}{3}n$    (or,  $ \gamma_{r5}(P(n,k)) =  \frac{5}{3}n$),
that have a 4RDF (5RDF) that is not obtained from a 3RDF of G ?  
Find (characterize) such graphs.

\medskip\noindent{\bf Question 2.}
Find (characterize) the cubic graphs for which 
  $ \gamma_{r5}(G) =  \frac{5}{3}n$ $\Rightarrow$  $ \gamma_{r4}(G) =  \frac{4}{3}n$ $\Rightarrow$ $ \gamma_{r3}(G) = n$.
Is there a counterexample ?

\medskip\noindent
Here we give an example that answers in part  the first  question.

\medskip\noindent{\bf Example.} 
Start the construction with 3 copies of an auxiliary graph $H$, a subvidived $K_4$. 
($H$ is obtained by replacing each edge of $K_4$ with  a path of length 2.)
Color the vertices of degree 3 of each $H$  with 4 different colors.
Add 6 vertices, and assign to them (all six) different two elemet subsets of $\{1,2,3,4\}$.
Connect each of the six new vertices to the vertices of degree two in each of  copy of $H$ (the subdivided $K_4$)
in a way that gives rise to a 4RDF. 
For example as follows: the vertex that is assigned $\{1,2 \}$ is connected to the vertex subdividing the edge between vertices 
colored by 3 and 4  (in each copy of $H$). This clearly defines a 4-rainbow domination assignement.
The graph constructed has $6+ 3\times 10 = 36$ vertices, 
and the number of colors used is $6\times 2 + 3\times 4 = 24$, hence   $\gamma_{r4} =  \frac{2}{3}|V(G)|$. 

It seems to be clear that this 4RDF can not be obtained from a  3RDF.
However, it is easy to find   a singleton 3RDF by suitably coloring the vertices that subdivide edges of $K_4$.
\qed 

\section{Conclusions}

The results  given in this paper in a  way conclude the  study of    rainbow domination of   generalized Petersen graphs 
$P(ck,k)$  because for all nontrivial $t$, we have characterizations of 
$P(ck,k)$ that attain general the lower bound, and we have good upper bounds for general case. 
In this paper we have provided the missing results for 4-rainbow and 5-rainbow domination
by proving Theorems   \ref{MainTheoremX2novi} and  \ref{MainTheoremX1}.

The upper bounds are derived from a relaxed versoin of the upper bounds for 3-rainbow domination. 
As indicated in the text, we could apply the same proof technique to elaborate on the more detailed upper bounds,
recalled in   Theorems \ref{MainTheoremOLD}, \ref{MainTheorem2OLD}, and \ref{MainTheorem3OLD}.
Straightforward proof, application of  Theorem    \ref{TheoremGeneralNOV} to each case separately,  would give :

 \begin{Theorem}  \label{MainTheoremDETAILS4}   
Bounds for  4-rainbow domination number $\gamma_{r4}$    of generalized Petersen graphs $P(ck,k)$   are: 
\begin{itemize}
\item For  $c \equiv 0\pmod 6$. 
 \begin{itemize}
 \item 
  If $ k \equiv 1,5\pmod 6$, then 
  $ \gamma_{r4}(P(ck,k))    =  \frac{4}{3}ck$;
 \item 
  If $ k \equiv 0\pmod 2$, then $
  \frac{4}{3}ck < \gamma_{r4}(P(ck,k))  \leq    \frac{5}{3}( c (k+ \frac{1}{2}) )$;
 \item 
  If $k \equiv 3\pmod 6$,  then  
 $\frac{4}{3}ck   < \gamma_{r4}(P(ck,k))    \leq   \frac{5}{3}( c(k+1) )$. 
 \end{itemize}

\item For  $c$ odd.  
 \begin{itemize}
 \item 
  { 
  If $ k \equiv 1,5\pmod 6$,} then 
  $ \frac{4}{3}ck <
  \gamma_{r4}(P(ck,k)) \leq   \frac{4}{3}(ck + \lceil\frac{k}{2} \rceil ) $;
 \item 
  { If $ k \equiv 0\pmod 2$,} then $\frac{4}{3}ck <\gamma_{r4}(P(ck,k))  \leq  \frac{4}{3}( ck+ \lfloor\frac{c}{2}\rfloor  + \frac{k}{2})$;
 \item 
  {  
  If $k \equiv 3\pmod 6$,}  then  
 $\frac{4}{3}ck   < \gamma_{r4}(P(ck,k))  \leq  \frac{4}{3}( c(k+1) + \lceil \frac{ k-2}{2} \rceil)$. 
 \end{itemize}

\item For  $c$  even, and $c \not\equiv 0\pmod 6$.
 \begin{itemize}
 \item 
  { 
  If $ k \equiv 1,5\pmod 6$,} then 
    $ \frac{4}{3}ck <   \gamma_{r4}(P(ck,k)) \leq \frac{4}{3}(ck + k + 1)$;
 \item 
  { If $ k \equiv 0\pmod 2$,} then $
  \frac{4}{3}ck < \gamma_{r4}(P(ck,k)) \leq  \frac{4}{3}(ck+ \frac{c}{2}  + k)$;
 \item 
  {  
  If $k \equiv 3~\pmod 6$,}  then  
 $\frac{4}{3}ck   < \gamma_{r4}(P(ck,k))  \leq  \frac{4}{3}( ck+  c + k-2)$. 
 \end{itemize} 
\end{itemize}
\end{Theorem}

 \begin{Theorem}  \label{MainTheoremDETAILS5}   
Bounds for  5-rainbow domination number $\gamma_{r5}$   of generalized Petersen graphs $P(ck,k)$   are: 
\begin{itemize}
\item For  $c \equiv 0\pmod 6$. 
 \begin{itemize}
 \item 
  If $ k \equiv 1,5\pmod 6$, then 
  $ \gamma_{r5}(P(ck,k))    =  \frac{5}{3}ck$;
 \item 
  If $ k \equiv 0\pmod 2$, then $
  \frac{5}{3}ck < \gamma_{r5}(P(ck,k))  \leq    \frac{5}{3}( c (k+ \frac{1}{2}) )$;
 \item 
  If $k \equiv 3\pmod 6$,  then  
 $\frac{5}{3}ck   < \gamma_{r5}(P(ck,k))    \leq  \frac{5}{3}( c(k+1))$. 
 \end{itemize}

\item For  $c$ odd.  
 \begin{itemize}
 \item 
  { 
  If $ k \equiv 1,5\pmod 6$,} then 
  $ \frac{5}{3}ck <
  \gamma_{r5}(P(ck,k)) \leq   \frac{5}{3}(ck + \lceil\frac{k}{2} \rceil ) $;
 \item 
  { If $ k \equiv 0\pmod 2$,} then $\frac{5}{3}ck <\gamma_{r5}(P(ck,k))  \leq  \frac{5}{3}( ck+ \lfloor\frac{c}{2}\rfloor  + \frac{k}{2})$;
 \item 
  {  
  If $k \equiv 3\pmod 6$,}  then  
 $\frac{5}{3}ck   < \gamma_{r5}(P(ck,k))  \leq  \frac{5}{3}( c(k+1) + \lceil \frac{ k-2}{2} \rceil)$. 
 \end{itemize}

\item For  $c$  even, and $c \not\equiv 0\pmod 6$.
 \begin{itemize}
 \item 
  { 
  If $ k \equiv 1,5\pmod 6$,} then 
    $ \frac{5}{3}ck <   \gamma_{r5}(P(ck,k)) \leq \frac{5}{3}(ck + k + 1)$;
 \item 
  { If $ k \equiv 0\pmod 2$,} then $
  \frac{5}{3}ck < \gamma_{r5}(P(ck,k)) \leq  \frac{5}{3}(ck+ \frac{c}{2}  + k)$;
 \item 
  {  
  If $k \equiv 3~\pmod 6$,}  then  
 $\frac{5}{3}ck   < \gamma_{r5}(P(ck,k))  \leq  \frac{5}{3}( ck+  c + k-2)$. 
 \end{itemize} 
\end{itemize}
\end{Theorem}  

It should be noted that the bounds given in the last two theorems are clearly not best possible. 
Namely, they are obtained from certain constructions (elaborated in \cite{Darja3RD}), and the method used assumes that 
one of the colors ir replaced by a set that includes the additional color (or, two colors). 
Thus it is very likely that improved bounds can be obtained, a question being whether there are elegant methods to do so. 

On the other hand, exact values, at least for some cases (e.g.  for $c$ small enough), seem to be tractable. 
For example using the suitably adapted algebraic method based on path algebras \cite{KLAVZAR199693}, 
see for example recent application of the method to a couple of related problems, namely 
the  independent rainbow domination on graphs $P(n,2)$ and $P(n,3) $  \cite{gabrovsek}
and the  2-rainbow independet domination of Cartesian products  of paths and cycles  \cite{gabrovsek2}.

\bigskip
\noindent{\bf Acknowledgement.}
This work was supported in part by the Slovenian Research Agency (grant  P2-0248).
 
%% References with BibTeX database:

\bibliographystyle{elsarticle-num}
%\bibliography{<your-bib-database>}
\bibliography{RainbowDominationJanezURLji}

%% Authors are advised to use a BibTeX database file for their reference list.
%% The provided style file elsarticle-num.bst formats references in the required Procedia style

%% For references without a BibTeX database:

% \begin{thebibliography}{00}

%% \bibitem must have the following form:
%%   \bibitem{key}...
%%

% \bibitem{}

% \end{thebibliography}

\end{document}